\documentstyle[12pt]{article}
\newcommand{\be}{\begin{equation}}
\newcommand{\ee}{\end{equation}}
\newcommand{\bea}{\begin{eqnarray}}
\newcommand{\eea}{\end{eqnarray}}
\newcommand{\nn}{\nonumber}
\newcommand{\nid}{\noindent}
\newcommand{\bx}{\mbox{\ \ \ $\Box$}}
\newcommand{\qe}{\mbox{exp${}_q$}}
\def\one{\raise0.0ex\hbox{$1$\kern-0.35em\raise0.0ex\hbox{$1$}}}
\def\C{\raise0.0ex\hbox{$\bf C$\kern-0.72em\raise0.0ex\hbox{$1$}}}
\def\R{\raise0.0ex\hbox{$\bf R$\kern-0.47em\raise0.0ex
\hbox{$\rule{0.8pt}{1.6ex}$}}}
\def\larrr{\raise0.3ex\hbox{$\longrightarrow$\kern-1.5em\raise-1.1ex
\hbox{$\scriptstyle{r\rightarrow 1}$}}}

\normalbaselines
\begin{document}
\begin{titlepage}

\begin{center}

\vspace*{1.0cm}

{\Large{\bf The Reversed q-Exponential Functional Relation.}}

\vskip 1.5cm

{\large {\bf David Fairlie}}\renewcommand{\thefootnote}
{\dagger}\footnote{E-mail address: david.fairlie@durham.ac.uk} \\
{\large {\bf Ming-Yuan Wu}}\renewcommand{\thefootnote}
{\ddagger}\footnote{E-mail address: vickiwu@ms5.hinet.net}

\vskip 0.5cm

Department of Mathematical Sciences, \\
University of Durham, \\
Durham DH1 3LE, \\
UK

\end{center}

\vspace{3 cm}

\begin{abstract}
\normalsize

After obtaining some useful identities, we  prove an additional
 functional relation for $q$ exponentials  with reversed order of 
multiplication, as well as the well known direct one in a completely 
rigorous manner. 

\end{abstract}
\end{titlepage}

\section{Introduction}

\nid One of the most appealing results to come out of $q$-analysis
 is that the $q$-exponential function, defined by 
\(
{}_qD_x \qe (x) = \qe (x) \, ,
\)
where ${}_qD_x$ is the $q$-derivative, also satisfies the same defining
functional relationship for ordinary exponential functions (up to
normalization), given by 
\be
F(x)F(y)=F(x+y) \, ,                                       \label{1}
\ee
provided that $xy=q^{-1}yx$. (that is, $(x,y)$ belongs to the Manin quantum
plane.) This result was first found by Sch\"utzenberger \cite{schu}) long
before the non-commutative aspects of q-analysis were generally recognised and
has been rediscovered many times subsequently e.g. in \cite{cigler,fairlie}.
It can be proved by means of $q$-combinatorics 
\cite{schu,cigler}, or by an argument based on the definition of the $q$-
exponential as an eigenfunction of the $q$-derivative \cite{fairlie}.

\nid Besides the above well-known result, there is, in fact, an additional
 functional relation in the opposite order for the $q$-exponential functions, 
 which is not so well known given by
\be
F(y)F(x)=F(\;x+y+(1-q^{-1})yx\;) \, ,                      \label{2}
\ee 
provided that the same condition $xy=q^{-1}yx$ holds. We first became aware
of this relationship in the work of  L. Faddeev and A. Yu Volkov in their 
study of lattice Virasoro algebra \cite{faddeev} who obtained a similar result
in the case of a different realisation of the $q$-exponential, in terms of an
infinite product. Their definition of the $q$-exponential
  suffered from the drawback that it did not go over into the ordinary 
exponential function in the commuting limit $q \rightarrow 1$.) 
In this paper, we will provide a completely rigorous proof  of the
reverse functional relation in the form stated in (\ref{2}).  
The proof 
is somewhat tricky in that a seemingly unrelated identity has to be 
obtained first as an intermediate step.   

\section{Proof of the Reversed q-Exponential Functional Relation}  

\nid For completeness we quickly review Sch\"utzenberger and Cigler's 
result, which will be used in our subsequent proof:
\be
\qe x\,\qe y=\qe (x+y)\qquad \mbox{if $xy=q^{-1}yx$}\, , \label{3}
\ee
where 
\[
\qe x\equiv \sum_{n=0}^{\infty}\frac{x^n}{[n]!}\qquad [n]\equiv 
\sum_{k=0}^{n-1}q^k\qquad [n]!\equiv [n][n-1]\cdots [1] \, .
\]
\vskip 0.9cm
\nid $\langle$\mbox{ Proof }$\rangle$
\vskip 0.1cm
\bea
\qe x\,\qe y & = & \left(\sum_{m=0}^{\infty}\frac{x^m}{[m]!}\right)\left(
                   \sum_{n=0}^{\infty}\frac{y^n}{[n]!}\right)\nn \\
             & = & \sum_{k=0}^{\infty}\sum_{r=0}^{k}\frac{x^ry^{k-r}}{[r]!
                   [k-r]!}\nn \\
             & = & \sum_{k=0}^{\infty}\frac{1}{[k]!}\left(\sum_{r=0}^{k}
                   \frac{[k]!}{[r]![k-r]!}x^ry^{k-r}\right)\nn \\
             & = & \sum_{k=0}^{\infty}\frac{(x+y)^k}{[k]!}\mbox{\ \ \ \  
(by (A1), see Appendix)}\nn \\
             & = & \qe (x+y)\nn \\
             &   & \mbox{Q.E.D.}\bx\nn
\eea
\vskip 0.6cm
\nid Now let us go on to prove the following formula,
\bea
x^n & = & 1+\sum_{r=1}^{n}\frac{[(q^{n-r+1}-1)(q^{n-r+2}-1)\cdots 
(q^n-1)][(x-1)(x-q)\cdots (x-q^{r-1})]}{(q-1)(q^2-1)\cdots (q^r-1)}\nn \\
  & \equiv & \sum_{r=0}^{n}\frac{[(q^{n-r+1}-1)(q^{n-r+2}-1)\cdots (q^n-1)]
              [(x-1)(x-q)\cdots (x-q^{r-1})]}{(q-1)(q^2-1)\cdots (q^r-1)}
               \, .                                      \label{4} 
\eea
\vskip 0.9cm
\nid $\langle$\mbox{ Proof }$\rangle$
\vskip 0.6cm
\nid Suppose for some $n=k$, we have 
\[
x^k=\sum_{r=0}^{k}\frac{[(q^{k-r+1}-1)(q^{k-r+2}-1)\cdots (q^k-1)][
          (x-1)(x-q)\cdots (x-q^{r-1})]}{(q-1)(q^2-1)\cdots (q^r-1)} \, .
\]
Now, consider $x^{k+1}$,
\bea
x^{k+1} & = & \sum_{r=0}^{k}\frac{[(q^{k-r+1}-1)(q^{k-r+2}-1)\cdots (q^k-1)][
              (x-1)(x-q)\cdots (x-q^r)]}{(q-1)(q^2-1)\cdots (q^r-1)}+\nn \\
        &   & \sum_{
              r=0}^{k}\frac{q^r[(q^{k-r+1}-1)(q^{k-r+2}-1)\cdots (q^k-1)][
              (x-1)(x-q)\cdots (x-q^{r-1})]}{(q-1)(q^2-1)\cdots (q^r-1)}\nn
              \\
        & = & (x-1)(x-q)\cdots (x-q^k)+\nn \\
        &   & \sum_{r=0}^{k-1}\frac{[(q^{k-r+1}-1)(q^{k-r+2}-1)\cdots 
              (q^k-1)]
              [(x-1)(x-q)\cdots (x-q^r)]}{(q-1)(q^2-1)\cdots (q^r-1)}+\nn \\
        &   & \sum_{r=0}^{k}\frac{q^r[(q^{k-r+1}-1)(q^{k-r+2}-1)\cdots 
              (q^k-1)]
              [(x-1)(x-q)\cdots (x-q^{r-1})]}{(q-1)(q^2-1)\cdots (q^r-1)}
              \nn 
\eea
\newpage
\bea
        & = & (x-1)(x-q)\cdots (x-q^k)+\nn \\
        &   & \sum_{r=1}^{k}\frac{[(q^{k-r+2}-1)(q^{k-r+3}-1)\cdots (q^k-1)]
              [(x-1)(x-q)\cdots (x-q^{r-1})]}{(q-1)(q^2-1)\cdots (q^{r-1}-1)}
              + \nn \\
        &   & \sum_{r=0}^{k}\frac{q^r[(q^{k-r+1}-1)(q^{k-r+2}-1)\cdots 
              (q^k-1)]
              [(x-1)(x-q)\cdots (x-q^{r-1})]}{(q-1)(q^2-1)\cdots (q^r-1)}
              \nn \\
        & = & \sum_{r=0}^{k+1}\frac{[(q^{(k+1)-r+1}-1)(q^{(k+1)-r+2}-1)\cdots 
             (q^{k+1}-1)]
             [(x-1)(x-q)\cdots (x-q^{r-1})]}{(q-1)(q^2-1)\cdots (q^r-1)}
              \, . \nn
\eea
Since for $n=1$, obviously we have 
\[
x=\sum_{r=0}^{1}\frac{[(q^{1-r+1}-1)(q^{1-r+2}-1)\cdots (q^1-1)]
  [(x-1)(x-q)\cdots (x-q^{r-1})]}{(q-1)(q^2-1)\cdots (q^r-1)} \, ,
\]
the proof is complete.\bx
\vskip 0.6cm
\nid There follows another identity which is a simple consequence of the
previous one; 
\be
\sum_{r=0}^{m\mbox{\tiny\,or\,}n}\frac{q^{\frac{r(r-1)}{2}-mn}(q-1)^r}{
[m-r]![n-r]![r]!}=\frac{1}{[m]![n]!}\, .                 \label{5}
\ee
\vskip 0.9cm
\nid $\langle$\mbox{ Proof }$\rangle$
\vskip 0.1cm
\bea
 &   &\sum_{r=0}^{n}\frac{q^{\frac{r(r-1)}{2}-mn}(q-1)^r}{[m-r]![n-r]![r]!}
      \nn \\
 & = &\sum_{r=0}^{n}\frac{q^{\frac{r(r-1)}{2}-mn}(q-1)^r
      ([m-r+1][m-r+2]\cdots [m])([n-r+1][n-r+2]\cdots [n])}{[m]![n]![r]!
      }\nn \\
 & = &\sum_{r=0}^{n}\frac{q^{\frac{r(r-1)}{2}-mn}\! [(q^{
      m-r+1}\! -\!\! 1)(q^{m-r+2}\! -\!\! 1)\cdots (q^m\! -\!\! 1)]
      [(q^{n-r+1}\! -\!\! 1)(q^{n-r+2}\! -\!\! 1)\cdots
      (q^n\! -\!\! 1)]}{[m]![n]!(q-1)(q^2-1)\cdots (q^r-1)}\nn \\
 & = &\frac{1}{[m]![n]!}\sum_{r=0}^{n}\frac{[(q^m-1)(q^m-q)\cdots 
      (q^m-q^{r-1})][(q^{n-r+1}-1)(q^{n-r+2}-1)\cdots (q^n-1)]}{(q^m)^n[
      (q-1)(q^2-1)\cdots (q^r-1)]}\nn \\
 & = &\frac{1}{[m]![n]!}\mbox{\ \ \ \ (by identity (4))} \, .\nn
\eea
\nid The proof is completed by noting that the above identity is 
symmetric in $m$ and $n$.\bx
\vskip 0.6cm
\nid Equipped with the above identity, we are now able to achieve the desired
result,
\be
\qe y\,\qe x=\qe [x+y+(1-q^{-1})yx]\qquad \mbox{if $xy=q^{-1}yx$}\, .
\label{6}
\ee
\vskip 0.9cm
\nid $\langle$\mbox{ Proof }$\rangle$
\vskip 0.1cm
\bea
 &   & \qe y\,\qe x\nn \\
 & = & \left(\sum_{m=0}^{\infty}\frac{y^m}{[m]!}\right)\left(
             \sum_{n=0}^{\infty}\frac{x^n}{[n]!}\right)\nn \\
 & = & \sum_{m=0}^{\infty}\sum_{n=0}^{\infty}y^mx^n\sum_{r=0}^{\min\{m,n\}}
       \frac{q^{\frac{r(r-1)}{2}-mn}(q-1)^r}{[m-r]![n-r]![r]!}
       \mbox{\ \ \ \ (by identity (5))}\nn \\
 & = & \sum_{m=0}^{\infty}\sum_{n=0}^{\infty}\sum_{r=0}^{\min\{m,n\}}
       \frac{q^{-r(n-r)}y^rx^{n-r}}{[n-r]!}\cdot
       \frac{q^{-\frac{r(r-1)}{2}}(1-q^{-1})^rx^r}{[r]!}\cdot
       \frac{y^{m-r}}{[m-r]!}\nn \\
 & = & \sum_{m=0}^{\infty}\sum_{n=0}^{\infty}\sum_{r=0}^{\min\{m,n\}}
       \frac{x^{n-r}}{[n-r]!}\cdot
       \frac{q^{-\frac{r(r-1)}{2}}(1-q^{-1})^ry^rx^r}{[r]!}\cdot
       \frac{y^{m-r}}{[m-r]!}\nn \\
 & = & \sum_{m=0}^{\infty}\sum_{n=0}^{\infty}\sum_{r=0}^{\min\{m,n\}}
       \frac{x^{n-r}}{[n-r]!}\cdot\frac{(1-q^{-1})^r(yx)^r}{[r]!}\cdot
       \frac{y^{m-r}}{[m-r]!}\nn \\
 & = & \left(\sum_{l=0}^{\infty}\frac{x^l}{[l]!}\right)\left(\sum_{
       k=0}^{\infty}\frac{[(1-q^{-1})yx]^k}{[k]!}\right)\left(\sum_{
       h=0}^{\infty}\frac{y^h}{[h]!}\right)\nn \\
 & = & \qe x\,\cdot \qe [(1-q^{-1})yx]\cdot \qe y \nn \\
 & = & \qe [x+(1-q^{-1})yx]\cdot \qe y\mbox{\ \ \ \ (by (3), as 
       $x(1-q^{-1})yx=q^{-1}(1-q^{-1})yxx$)}\nn \\
 & = & \qe [x+(1-q^{-1})yx+y]\nn \\
 &   & \mbox{\hspace{4cm}(by (3), as $[x+(1-q^{-1})yx]y=q^
       {-1}y[x+(1-q^{-1})yx]$)}\nn \\
 & = & \qe [x+y+(1-q^{-1})yx]\nn \\
 &   & \mbox{Q.E.D.}\bx\nn
\eea

\vspace{0.3cm}

\nid{\Large{\bf Appendix}}

\vspace{0.3cm}

\nid The following is the so-called $q$-binomial expansion formula and its
proof: 
\[
\mbox{\hspace{4.5cm}}
(x+y)^n=\sum_{r=0}^{n}\left[\!\! 
\begin{array}{c}
n\\
r
\end{array}\!\!\right]x^ry^{n-r} \, ,\mbox{\hspace{4.5cm}(A1)}
\]
where
\[
\left[\!\!
\begin{array}{c}
n\\
r
\end{array}\!\!\right] \equiv \frac{[n]!}{[r]![n-r]!} \qquad 
[n] \equiv \sum_{k=0}^{n-1}q^k \qquad [0]! \equiv 1 \, ,
\]
subject to the condition that $xy=q^{-1}yx$, $q$ being some complex number. 
\vskip 0.9cm
\nid $\langle$\mbox{ Proof }$\rangle$ 
\vskip 0.6cm
\nid Suppose for some $n=k$, we have
\[
(x+y)^k=\sum_{r=0}^{k}\left[\!\!
\begin{array}{c}
k\\
r
\end{array}\!\!\right]x^ry^{k-r} \, .
\] 
Now consider $(x+y)^{k+1}$,
\bea
(x+y)^{k+1} & = & \sum_{r=0}^{k}\frac{[k]!}{[r]![k-r]!}(x+y)x^ry^{k-r}\nn \\
            & = & \sum_{r=0}^{k}\frac{[k]!}{[r]![k-r]!}x^{r+1}y^{k-r}
                  +\sum_{r=0}^{k}\frac{[k]!}{[r]![k-r]!}q^rx^ry^{k-r+1}\nn \\
            & = & x^{k+1}+\sum_{r=0}^{k-1}\frac{[k]!}{[r]![k-r]!}x^{r+1}
                  y^{k-r}+\sum_{r=0}^{k}\frac{[k]!}{[r]![k-r]!}q^rx^r
                  y^{k-r+1}\nn \\
            & = & x^{k+1}+\sum_{r=1}^{k}\frac{[k]!}{[r-1]![k-r+1]!}x^ry^{k-r
                  +1}+\sum_{r=0}^{k}\frac{[k]!}{[r]![k-r]!}q^rx^ry^{k-r+1}
                  \nn \\
            & = & x^{k+1}+\sum_{r=1}^{k}\frac{[k]!(1+q+\cdots +q^{r-1})}{[r]!
                  [k-r]!(1+q+\cdots +q^{k-r})}x^ry^{k-r+1}+ \nn \\
            &   & \sum_{r=0}^{k}\frac{[k]!}{[r]![k-r]!}q^rx^ry^{k-r+1}\nn \\
            & = & x^{k+1}+\sum_{r=0}^{k}\frac{[k]!(1+q+\cdots +q^k)}{[r]!
                  [k-r]!(1+q+\cdots +q^{k-r})}x^ry^{k-r+1}\nn \\
            & = & \sum_{r=0}^{k+1}\frac{[k+1]!}{[r]![k+1-r]!}x^ry^{k+1-r}\,
                  , \nn
\eea
so the same formula holds for $n=k+1$. \\
Since for $n=1$, obviously we have 
$\displaystyle x+y=\sum_{r=0}^{1}\frac{[1]!}{[r]![1-r]!}x^ry^{1-r}$,
the proof is complete.\bx
    
\vspace{0.35cm}

\nid {\large{\bf Acknowledgements}}

\vspace{0.2cm}

\nid M Y Wu is grateful to members in the Department of Mathematical Sciences 
of the University of  Durham for all relevant help.     This work was
partially supported by a British ORS Award as well as a Durham University 
Research Award.

\end{document}